\input amstex
\documentstyle{ams-j}

\define\e{\varepsilon}
\redefine\o{\omega}
\redefine\L{\Lambda}
\define\vp{\varphi}
\define\Hom{\operatorname{Hom}}
\define\Ind{\operatorname{Ind}}
\define\supp{\operatorname{supp}}

\NoBlackBoxes

\topmatter
\title{On isotropic Bessel models for principal series representations of $SO_{n,n+1}(\Bbb
R)$}\endtitle
\author
David Goldberg
\endauthor
\address Department of Mathematics, Purdue University, West Lafayette, IN 47907
\endaddress
\email goldberg\@math.purdue.edu\endemail
\abstract
We show that principal series representations for $SO_{n+1,n}(\Bbb R)$ have non-vanishing
rank one Bessel models of a very particular type.  The Bessel functionals are explicitly constructed with respect to an $SO_{2,1}(\Bbb R)$--isotropic vector in $\Bbb R^3.$
\endabstract
\subjclass Primary 22E500,  Secondary 22E35\endsubjclass
\endtopmatter

\document

\subheading{Introduction}
Whittaker models have been crucial in the development of the theory of automorphic
forms.  Within the Langlands-Shahidi method, Whittaker models give rise to the definition of 
local coefficients \cite{\bf Sh1,Sh2}, which proved critical in establishing
 the proof of Langlands's conjecture relating local $L$--functions
to Plancherel measures for generic representations \cite{\bf Sh3}.  This, in turn, led to the
recent breakthroughs in functoriality \cite{\bf K,KS}, the transfer of generic cusp
forms from classical groups to general linear groups \cite{\bf CKPSS1,CKPSS2}, and progress on the local Langlands conjecture for $SO_{2n+1}$ \cite{\bf JS}.

This raises the question of whether other models can be used to obtain similar results for the 
non-generic spectrum.  In \cite{\bf FG} we used the generalized Bessel models of
\cite{\bf GPSR} to give some preliminary results on local coefficients for non-generic
representations of certain classical groups.  There, the analog of Rodier's Theorem is
established, subject to a minimality condition, and this leads to results on local coefficients of
the type in \cite{\bf Sh1}.  In order to carry the philosophy of the Langlands-Shahidi method
forward for non-generic representations of these groups, several questions need to be answered.  In particular, the minimality
condition has to be relaxed or removed,  the analog of the results of \cite{\bf Sh2} need to
be established for  these models, and global methods will need to be exploited to obtain results
relating Plancherel measures and $L$--functions.

In addressing the questions on local coefficients for Bessel models of groups over archimedean
fields, several issues arise, including the existence and uniqueness of such models.  Here we
address the existence of Bessel models for principal series of $G=SO_{n,n+1}(\Bbb R).$  By the Casselman Embedding Theorem, every irreducible representation is a subquotient of such a principal series.  Thus, we wish to investigate to what extent principal series have generalized Bessel models. 
We show that every such representation has a Bessel model with respect to a very
particular character of a subgroup of $SO_{1,2}(\Bbb R).$    In the language  \cite{\bf FG}
this model is rank one, and as remarked there \cite{\bf loc cit, Remark 2.10}, the minimality condition
is automatically satisfied.  Thus, one can hope to carry out the steps outlined above.  The
Rodier result can be establsihed for these models, i.e., one can establish that the dimension of the 
space of Bessel models with compatible parameters is preserved under parabolic induction. 
We also hope to prove the existence of local
coefficient in this case, and to show that these local coefficients are given by Artin
$L$--functions. 
However, it seems unlikely that the models described in this note will be sufficient to handle the most general case.
We expect that models of all types and ranks may be necessary to fully generalize the work of \cite{\bf Sh2}.
All this will be the subject of further work.

The author wishes to thank Jeff Adams, Sol Firedberg, Freydoon Shahidi, and Ramin Takloo-Bighash for discussions useful in this work.

 \subheading{Section 1: Preliminaries}

Let $G=SO_{n,n+1}(\Bbb R),$ defined with respect to the form
$$J_{2n+1}=\pmatrix &&&&1\\&&&.\\&&.\\&.\\1\endpmatrix.$$ 
 We Let $B=M_0A_0U$ be the Borel subgroup
of upper triangular matrices in $G,$  and its Langlands decomposition.
Then,
$$M_0=\left\{\pmatrix
\delta_1\\&\delta_2\\&&\ddots\\&&&\delta_n\\&&&&1\\&&&&&\delta_n^{-1}\\&&&&&&\ddots\\&&&&&&&\delta_1^{-1}\endpmatrix\Bigg|\,
\delta_i\in\left\{\pm 1\right\}\right\}\simeq\Bbb Z_2^n$$
and 
$$A_0=\left\{\pmatrix x_1\\&x_2\\&&\ddots\\&&&
x_n\\&&&&1\\&&&&&x_n^{-1}\\&&&&&&\ddots\\&&&&&&&x_1^{-1}\endpmatrix\Bigg|\, x_i\in\Bbb
R^+\right\}.$$
The unipotent radical $U$ is the set of upper triangular unipotent matrices in $G.$
We denote by $\Phi$ the
roots of
$A_0$ in $G,$ and by $\Delta=\{e_1-e_2,\dots e_{n-1}-e_n,e_n\}$ the simple roots. Let $\rho$
be half the sum of the positive roots.  We denote by
$W=W(G,A_0)=N_G(A_0)/A_0$  the Weyl group.

Now let $U_1\subset U$ be the group of unipotent upper triangular
matrices whose middle $3\times 3$ block is the indentity matrix.  
Define 
$$M_1=Z_G(U_1)=\left\{\pmatrix I_{n-1}\\&h\\&&I_{n-1}\endpmatrix\big|h\in
SO_{1,2}(\Bbb R)\right\}.$$
Let
$\psi:\Bbb R\rightarrow\Bbb C^\times$ be given by $\psi(x)=e^{2\pi i x}.$
Fix a non-zero element $\bold a=(a,b,c)$ of $\Bbb R^3$ and define a
character $\chi_{\bold a}=\chi$ of $U_1$ by
$$\chi((u_{ij})=\psi(u_{12}+u_{23}+\dots +u_{n-2,n-1}+au_{n-1,n}+bu_{n-1,n+1}+cu_{n-1,n+2}).$$
Let $M_{\chi}=\{g\in M_1|\chi^g=\chi\}.$ Then $M_\chi$ is isomorphic to the
stabilizer of $\bold a$ in $SO_{1,2}(\Bbb R).$
Set $R_\chi=M_\chi U_1.$   For an irreducible (Frechet)
representation $\o$ of $M_\chi,$ we denote by $\o_\chi$ its extension to
$R_\chi$ by $\chi.$  

\definition{Definition 1.1} We say that a Frechet representation $\Pi$ of $G$ has an
$\o_\chi$--Bessel model if $\Hom_G(\Pi,\Ind_{R_\chi}^G(\o_\chi))\neq 0.$ 
Here, $\Ind$ means diffrentiably induced. 
\enddefinition

Let $\frak a^*$ be the  dual of the Lie algebra $\frak a$ of $A_0,$ and set
$\frak a^*_{\Bbb C}=\frak a^*\otimes_{\Bbb R}\Bbb C$ to be the complexification of $\frak a^*.$
If $\eta\in \hat M_0$ and $\nu\in\frak a^*_\Bbb C,$ we let $\pi_{\nu,\eta}=\Ind_B^H(\eta e^{\nu}).$  Let
$V(\nu,\eta)$ be the space of $\pi_{\nu,\eta}.$    Thus,
$$V(\nu,\eta)=\left\{ f\in C_c^\infty(G)\, |f(m_0a_0ug)=\eta(m_0)^{-1}e^{-(\nu+\rho)\log
a_0}f(g)\text{ for all } m_0a_0u\in B, g\in G\right\}.$$  We will show that these
representations  always have Bessel models with respect to a particular character $\chi.$

\subheading{Section 2 The main theorem}

We consider the special case where
$\chi=\chi_v$ with $v=\pmatrix 1\\-1\\-1/2\endpmatrix.$
Note that $v$ is isotropic with respect to our form.  If $v_0=\pmatrix 1\\0\\0\endpmatrix,$
then we can choose $g_0\in SO_{1,2}(\Bbb R)$ with $g_0v_0=v.$ In particular, we may take
$$g_0=\pmatrix1&0&0\\-1&1&0\\-1/2&1&1\endpmatrix.$$
Note that the stablizier $M_{v_0}$ of $v_0$ is  given by
$$M_{v_0}=\left\{\pmatrix
I_{n-1}\\&1&x&-x^2/2\\&0&1&-x\\&0&0&1\\&&&&I_{n-1}\endpmatrix\bigg| x\in\Bbb
R\right\}=U_{e_n},$$
which is the unipotent subgroup associated to the short simple root $e_n.$ Then
$M_{\chi}=M_v=g_0M_{v_0}g_0^{-1},$ and we find by direct computation that
$$M_\chi=\left\{\pmatrix
I_{n-1}\\&(\frac
x2+1)^2&\frac{x^2}{2}+x&-\frac{x^2}2\\&-\frac{x}2(1+\frac{x}2)&1-\frac{x^2}2&\frac{x^2}2-x\\
&-\frac{x^2}8&\frac{x}2(1-\frac{x}{2})&(\frac x2-1)^2\\&&&&I_{n-1}\endpmatrix\bigg|\, x\in\Bbb
R\right\}.$$  For $x\in\Bbb R,$ we let 
$$m(x)=\pmatrix
I_{n-1}\\&(\frac
x2-1)^2&\frac{x^2}{2}+x&-\frac{x^2}2\\&-\frac{x}2(1+\frac{x}2)&1-\frac{x^2}2&\frac{x^2}2-x\\
&-\frac{x^2}8&\frac{x}2(1-\frac{x}{2})&(\frac x2-1)^2\\&&&&I_{n-1}\endpmatrix,$$
and let  $X_{e_n}(x)=g_0^{-1}m(x)g_0.$  
As $M_\chi$ is abelian, any irreducible representation of $M_\chi$ is  a character.
We fix such a character $\omega,$ and note that it is of the form $\o(m(x))=e^{2\pi icx}$ for
some $c\in[0,1).$

\definition {Definition 2.1}  Let $\chi$ and $\omega$ be as above.  Define, formally, a
functional $\L_{\chi,\o}$ on $V(\nu,\eta)$ by
$$\L_{\chi,\o}(f)=\int_{R_\chi}f(w_0^{-1}mu)\chi(u)^{-1}\o(m)^{-1}\, dm\, du,\tag 2.1$$
where
$w_0=\pmatrix &&J_n\\&-1\\J_n\endpmatrix$ 
represents the longest element in the Weyl group.
\enddefinition

Let $\delta_{R_\chi}$ be the modular function on $R_\chi.$
We note that if $u_1\in U_1,$ then
$$\gather
\L_{\chi,\o}(\pi_{\nu,\eta}(u_1)f)=\int_{M_{\chi}}\int_{U_1}f(w_0^{-1}muu_1)\o(m^{-1})\chi(u)^{-1}\,du\,dm\\
=\int_{M_{\chi}}\int_{U_1}f(w_0^{-1}mu)\chi(uu_1^{-1})^{-1}\o(m^{-1})\,du\,dm=\chi(u_1)\L_{\chi,\o}(f).
\endgather$$
Also, if $m_1\in M_{\chi}$ then
$$\gather
\L_{\chi,\o}(\pi_{\nu,\eta}(m_1)f)=\int_{M_{\chi}}\int_{U_1}f(w_0^{-1}mum_1)\o(m^{-1})\chi(u)^{-1}\,du\,dm\\
=\int_{M_{\chi}}\int_{U_1}f(w_0^{-1}mm_1(m_1^{-1}um_1))\o(m^{-1})\chi(u)^{-1}\,du\,dm\\
=\delta_{R_{\chi}}(m_1)^{1/2}\int_{M_{\chi}}\int_{U_1}f(w_0^{-1}mm_1u)\o(m^{-1})\chi^{m_1}(u)^{-1}\,du\,dm\\
=\delta_{R_{\chi}}(m_1)^{1/2}\int_{M_{\chi}}\int_{U_1}f(w_0^{-1}mu)\o((mm_1^{-1})^{-1})\chi(u)^{-1}\,du\,dm\\
=\delta_{R_{\chi}}(m_1)^{1/2}\o(m_1)\L_{\chi,\o}(f).
\endgather$$
Now define $B_{\chi,\eta,f}(g)=\L_{\chi,\o}(\pi_{\nu,\eta}(g)f).$  Then the map $f\mapsto
B_{\chi,\eta,f}$ is (formally) a Bessel model for $\pi_{\nu,\eta}.$

\proclaim{Proposition 2.2}   The integral (2.1) converges absolutely for all $\nu$ in a half space, with meromorphic 
continuation to $\frak a_\Bbb C^*.$ This, then defines
a Bessel functional on $V(\nu,\eta).$
\endproclaim

\demo{Proof}
 
Note that
$$|\L_{\chi,\o}(f)|\leq\int_{M_\chi U_1}|f(w_0^{-1}mu)|\, dm\, du.\tag 2.2$$
We note that $g_0$ normalizes $U_1.$
Then we rewrite the right hand side of (2.2) as
$$\gather
\int_{U_{e_n}}\int_{U_1}|f(w_0^{-1}g_0X_{e_n}(x)g_0^{-1}u)|\, du\,
dx=\int_{U_{e_n}}\int_{U_1}|f(w_0^{-1}g_0X_{e_n}(x)ug_0^{-1})|\,du\,dx\\
=\int_{U}|(\pi_{\nu,\eta}(g_0)^{-1}f(X_{e_n}(-1)w_0^{-1}ug_0)|\, du=\int_U|f(w_0^{-1}ug_0^{-1})|\, du.\tag
2.3
\endgather$$

Now note that, if $A(\nu,\eta, w_0)$ is the standard intertwining operator, then 
$$A(\nu,\eta,w_0)f(g_0^{-1})=\int_U f(w_0^{-1}ug_0^{-1})\, du,$$ converges absolutely for $\nu$ satisfying 
$\operatorname{Re}\left<\nu, \alpha\right>\, >0$ for all $\alpha>0,$
with meromorphic continuation to all $\frak a_{\Bbb C}^*.$  Now (2.3) 
clearly also converges absoutely for this same half space, and continues meromorphically to all of $\frak
a_{\Bbb C}^*.$
\qed
\enddemo

We need to show that $\L_{\chi,\omega}\not\equiv 0.$
Let $\vp\in C_c^\infty(H)$ and define
$$f_\vp(g)=\sum_{M_0}\int_{A_0U}\eta(m_0)^{-1}e^{-(\nu+\rho)\log
a_0}\vp(m_0^{-1}a_0^{-1}u_0^{-1}g)\,du_0\,da_0.$$
Then $f_\vp\in V(\nu,\eta).$   Note that
$$\gather
\L_{\chi,\o}(f_\vp)=\int_{M_\chi}\int_{U_1}f_\vp(w_0^{-1}mu_1)\chi(u_1^{-1})\o(m^{-1})\,du_1\,dm\\
=\sum_{M_0}\eta(m_0)^{-1}\int_{M_{\chi}}\int_{U_1}\int_{A_0U}e^{-(\nu+\rho)\log
a_0}\vp(m_0^{-1}a_0^{-1}u_0^{-1}w_0^{-1}mu_1)\o(m)^{-1}\chi(u_1)^{-1}\,duda_0du_1dm.
\endgather$$

Let $P=MU_1$ be the standard parabolic subgroup with Levi component
$$M=\left\{\pmatrix x\\&g\\&&\varepsilon(x)\endpmatrix\,\Bigg| x\in GL_{n-1}(\Bbb R),g\in
SO_{1,2}(\Bbb R),\right\}$$ and unipotent radical $U_1.$ Then $M_1\subset M,$ so $M_\chi
U_1\subset P.$  Here $\e(x)=w_{n-1}\,^tx^{-1} w_{n-1}^{-1}.$
Note that $Bw_0B\subset Bw_0P$ is the open Bruhat cell in the $B\backslash G/P$--decomposition
(and $w_0=w_0^{-1}$).  Thus, $M_0A_0Uw_0^{-1}M_\chi U_1\subset Bw_0P.$ So, if
$\supp\vp\cap Bw_0P=\emptyset,$ then $\L_{\chi,\o}(f_\vp)=0.$

\proclaim{Lemma 2.2}  The map $B\times M_\chi U_1\rightarrow Bw_0M_\chi U_1$ given by
$(b,mu_1)\mapsto bw_0u_1$ is an isomorphism.
\endproclaim

\demo{Proof}  Suppose $b_1w_0m(x_1)u_1=b_2w_0m(x_2)u_2.$  Then,
$w_0^{-1}b_2^{-1}b_1w_0=m(x_2)u_2u_1^{-1}m(x_1)^{-1}\in \bar B\cap M_\chi U_1,$ where $\bar B$ is
the Borel opposite to $B.$  But noting the block form of the product, and since $M_\chi\cap\bar B=\{1\},$
we have $\bar B\cap M_\chi U_1$ is also trivial, which proves the
lemma.\qed\enddemo

Recall that $R_\chi=M_\chi U_1.$  Let $\vp\in C_c^\infty(Bw_0R_\chi)$	be of the form
$\vp(bw_0r)=\vp_B(b)\vp(w_0)\vp_{R_{\chi}}(r),$ with $\vp_B\in C_c^\infty(B),$ and
$\vp_{R_\chi}\in C_c^\infty(R_\chi).$ We further assume that
$\vp_{R_\chi}(m(x)u_1)=\vp_M(x)\vp_{U_1}(u_1),$ with $\vp_M\in C_c^{\infty}(\Bbb R),$ and
$\vp_{U_1}\in C_c^\infty(U_1).$

For $b=m_0a_0u_0\in M_0A_0U,$ we let $\eta(b)=\eta(m_0)$ and $\log b=\log a_0.$
Then 
$$\gather
\L_{\chi,\o}(f_\vp)=\int_B\int_{\Bbb R}\int_{U_1}\eta(b)^{-1}e^{-(\nu+\rho)\log
b}\vp_B(b^{-1})\vp(w_0)\vp_M(x)\o(x)^{-1}\vp_{U_1}(u_1)\chi(u_1)^{-1}du_1\,dx\,db\\
=\left(\int_B\vp_B(b^{-1})\eta(b)^{-1}e^{-(\nu+\rho)\log
b}\,db\right)\vp(w_0)\left(\int_{\Bbb R}\vp_M(x)
\o(x)^{-1}\,dx\right)\int_{U_1}\vp_{U_1}(u_1)\chi(u_1)^{-1}\,du_1\\
=\vp(w_0^{-1})\left(\int_B\vp_B(b^{-1})\eta(b)^{-1}e^{-(\nu+\rho)\log
b}\,db\right)\hat\vp_{M}(\o)\hat\vp_{U_1}(\chi),
\endgather$$
where $\hat\vp_{M}$ and $\hat\vp_{U_1}$ are Fourier transforms.

We can choose $\vp_M$ and $\vp_{U_1}$ so that their Fourier transforms at $\o$ and $\chi,$
respectively are non-zero.  Also, choose $\vp_B$ to have small enough (compact) support so that
the first factor is non-zero.  Then $\L_{\chi,\o}(f_\vp)\neq0,$ so we have a non-zero Bessel
functional for $\pi_{\nu,\eta}.$  

We now state our full result.

\proclaim{Theorem 2.3}  For every $\eta$ and $\nu$ there is an $\o\chi$--Bessel model for
$\pi_{\nu,eta},$ with
$\chi=\chi_v,$ and $v=\pmatrix 1\\-1\\-1/2\endpmatrix.$\qed
\endproclaim

\enddocument
\end